\documentclass[10pt]{article}
\usepackage[english]{babel}
\usepackage{amssymb,amsmath,amsthm,amsfonts}
\usepackage{graphicx}

\begin{document}

\noindent {ISSN 1066-369X, Russian Mathematics, 2018, Vol. 62, No.
2, pp. 28-33. }

\bigskip
\hrule height 1pt\kern 2pt \hrule height 0.25pt

\vspace{15ex}

\begin{center}
\large{\textbf{Vector Hamiltonians in Nambu mechanics}}

\bigskip

V.N.Dumachev\footnote{Voronezh institute of the MVD of Russia.
e-mail: dumv@comch.ru}
\end{center}

\begin{abstract}
We give a generalization of the Nambu mechanics based on vector
Hamiltonians theory. It is shown that any divergence-free phase flow
in $\mathbb{R}^n$ can be represented as a generalized Nambu
mechanics with $n-1$ integral invariants. For the case when the
phase flow in $\mathbb{R}^n$ has $n-3$ or less first integrals, we
introduce the Cartan concept of mechanics. As an example we give the
fifth integral invariant of Euler top.
\end{abstract}

\textbf{DOI:} 10.3103/S1066369X18020044


\section*{Introduction}

Consider phase flow of divergence-free type
\[
 \overset{\cdot}{\textbf{x}}=\textbf{l(x)}.\eqno{(1)}
\]
In this paper, all objects are supposed to be smooth. We denote by
$\Lambda^k(\mathbb{R}^n)$the space of $k$-differential forms on
$\mathbb{R}^n$, by $\Omega \in \Lambda^n(\mathbb{R}^n)$ a volume
form, by $T^k(\mathbb{R}^n)$ the space of k-vector fields, by
$H^k(\mathbb{R}^n)$ kth group of de Rham cohomologies, by $[\omega]$
the class of cohomologies of the form $\omega$. Henceforth we
suppose $\dim H^{n-1}(\mathbb{R}^n)=0$. Then $[\omega] \in
H^k(\mathbb{R}^n)$ means that $\omega$ is closed, and $[\omega]=0$
implies $\omega=d\nu$. Inner product of $X \in T(\mathbb{R}^n)$ and
$\omega \in \Lambda^k(\mathbb{R}^n)$ is $X \rfloor \omega \in
\Lambda^{k-1}(\mathbb{R}^n)$.  Vector field $X_\textbf{I}$ on
symplectic manifold $(\mathbb{R}^2,\Omega)$ is said to be
Hamiltonian if the corresponding 1-form $\Theta =X_\textbf{I}\rfloor
\Omega$ is closed. Due to Liouville's theorem, any Hamiltonian field
preserves the form of phase volume $\Omega \in
\Lambda^n(\mathbb{R}^n)$, i.e., the Lie derivative of the form
$\Omega$ in vector field $X_\textbf{I}$ is zero:
\[
L_{X_\textbf{I}}\Omega=X_\textbf{I}\rfloor d\Omega+d(X_\textbf{I}
\rfloor \Omega)=0.
\]
In the classic case, $\Omega \in \Lambda^2(\mathbb{R}^2)$, i.e.,
$d\Omega=0$, and provided the condition of Poincar\'e's lemma, from
$d(X_\textbf{I} \rfloor \Omega)=0$ follows $X_\textbf{I} \rfloor
\Omega=d\textbf{I} \in \Lambda^1(\mathbb{R}^2)$, i.e., the function
$\textbf{I} \in \Lambda^0(\mathbb{R}^2)$ is an invariant of
dynamical system $X_\textbf{I}$. The phase flow (1) can be
represented by means of Poisson bracket in the form
\[
 \overset{\cdot }{\textbf{x}}=\{\textbf{I},\textbf{x}\}.
\]

\section{Nambu phase flows}

In generalization of this case, Y. Nambu \cite{1} supposes that
$n$-dimensional phase flow is described by means of $n-1$ invariants
\[
 \overset{\cdot }{\textbf{x}}=\{\textbf{I}_1,\textbf{I}_2,...,\textbf{I}_{n-1},\textbf{x}\}.\eqno{(2)}
\]
From geometric point of view, the solution to system (2) is an
integral curve $l=l(\mathbb{R}^n)$, i.e., an object of dimension
$\dim(l)=1$. A curve is said to be algebraic, if it can be realized
as complete intersection of $n-1$ hypersurfaces, i.e., $l= \bigcap
\limits_{k=1}^{n-1} \textbf{I}_k$. A hypersurface in algebraic
geometry is an object with dimensionless by one than that of
surrounding space, thus, any polynomial in $\mathbb{R}^n$ forms an
algebraic hypersurface \cite{2}. If invariant $\textbf{I}$ has a
physical sense (e.g., energy integral), it is called Hamiltonian
$\textbf{H}$. For instance, for energy integral
$\textbf{H}=\frac{1}{2}(x_1^2+x_2^2),$ from the preservation law

\noindent \begin{minipage}[t]{85mm}
\[
\overset{\cdot}{\textbf{H}}=0, \qquad \text{i.e.,} \qquad
\textbf{H}_{x_1}\overset{\cdot
}{\textbf{x}}_1+\textbf{H}_{x_2}\overset{\cdot}{\textbf{x}}_2=0
\]
the Hamilton motion equations follow
\[\overset{\cdot
}{\textbf{x}}_1=\textbf{H}_{x_2}, \; \overset{\cdot
}{\textbf{x}}_2=-\textbf{H}_{x_1},\]
which coincide with (1) in
$\mathbb{R}^2$. Hence, invariants of Eq.(1) are hypersurfaces, whose
complete intersection gives integral curve $l$. For instance, Nambu
phase flow in $\mathbb{R}^3$ is formed as tangent one to two
invariants $(\textbf{I}_1,\textbf{I}_2)$, i.e.,
\end{minipage} \hfill
\begin{minipage}[t]{70mm}
{\;}

\begin{picture}(60,130)

\qbezier(55,22)(90,70)(160,97)

\put(100,66){\vector(0,1){23}}

\put(100,66){\vector(3,-1){20}}

\qbezier(0,40)(100,60)(200,0)

\qbezier(0,40)(50,80)(100,100)

\qbezier(100,100)(150,80)(200,0)

\qbezier(0,80)(80,60)(120,-10)

\qbezier(0,80)(40,120)(80,140)

\qbezier(80,140)(160,100)(120,-10)

\put(150,40){$\textbf{I}_1$}

\put(85,80){$\textbf{N}_1$}

\put(50,100){$\textbf{I}_2$}

\put(110,50){$\textbf{N}_2$}

\put(130,100){$\textbf{l}=\textbf{N}_1 \times \textbf{N}_2$}

\end{picture}

\end{minipage}

\[ \overset{\cdot
}{\textbf{x}}=\textbf{l}=\textbf{N}_1 \times \textbf{N}_2,\] where
$\textbf{N}=\text{grad}\,\textbf{I}$. From the point of view of
exterior differential algebra $\textbf{I} \in
\Lambda^0(\mathbb{R}^3)$, $d\textbf{I}=(\textbf{N} \cdot
d\textbf{r}) \in \Lambda^1(\mathbb{R}^3)$, $\omega=d\textbf{I}_1
\wedge d\textbf{I}_2 \in \Lambda^2(\mathbb{R}^3)$. Denote
$\omega=(\textbf{l} \cdot d\textbf{S})$, and notice that $\omega$ is
always closed. It is straightforward that $\omega=\textbf{l}\rfloor
\Omega$, and then $d(\textbf{l}\rfloor \Omega)=0$,  and hence,
$\text{div}\, \textbf{l}=0$. Increase in dimension leads to the
following

\bigskip
\textbf{Theorem 1.} Nambu phase flow is divergence-free.

\bigskip
\textbf{Proof.} Consider differential form $\omega=\bigwedge
\limits_{k=1}^{n-1}d\textbf{I}_k$. It is straightforward to see that
it can be represented in the form $\omega=(
\overset{\cdot}{\textbf{x}} \cdot d\textbf{S})$, where $dS_k=(-)^k
dx_1 \wedge dx_2 \wedge ...  \wedge [dx_k] \wedge ... \wedge dx_n$
are Pl\"ucker coordinates of elementary hyperplatform, spanned on
increment vectors $(dx_1,...dx_n)$. In similar case further, we will
say that we represent system (1) as differential form $\omega$.
Formula $\omega=\bigwedge \limits_{k=1}^{n-1}d\textbf{I}_k$ implies
that $\omega$ is closed and $\text{div}\,
\overset{\cdot}{\textbf{x}}=0.$ $\Box$

Hence due to Poincar\'e lemma $\omega=d\textbf{h}, \; \textbf{h} \in
\Lambda^{n-2}(\mathbb{R}^n)$, and for Nambu phase flow
\[
d\textbf{h}=d\textbf{I}_1 \wedge ... \wedge d\textbf{I}_{n-1}.
\]
The quantity $\textbf{h}$ is called vector (tensor) Hamiltonian
\cite{3}. Further we will study system of the form
\[
d\textbf{h}=\bigwedge \limits \textbf{J}^k, \qquad \sum k =n-1,
\eqno{(3)}
\]
where $\textbf{J}^k \in \Lambda^k(\mathbb{R}^n)$.

In mechanics, quantity $\textbf{I} \in \Lambda^0(\mathbb{R}^n)$ is
called the first integral, $\textbf{J} \in \Lambda^1(\mathbb{R}^n)$
integral invariant, $\textbf{J} \in \Lambda^k(\mathbb{R}^n)$
generalized integral invariant \cite{4}.

\newpage
\section{Splitting de Rham cohomologies}\label{sec3}

Endow de Rham complex
\[
0\rightarrow \Lambda ^{0}(\mathbb{R}^n)\rightarrow \Lambda
^{1}(\mathbb{R}^n )\rightarrow ...\rightarrow \Lambda
^{n-1}(\mathbb{R}^n)\rightarrow \Lambda
^{n}(\mathbb{R}^n)\rightarrow 0
\]
with differential module $\{C,d\}$, then
\[
Z(C)=\ker d= \{x \in C|\; dx=0 \}\] are co-cycles of module
$\{C,d\}$ (the space of closed forms),
\[B(C)={\rm im}\, d= dC=\{x=dy|\; y \in C\}\]
are co-boundaries of module $\{C,d\}$ ((the space of exact forms).
In these notations, group of ith cohomologies of $H^i$ is a quotient
of ith co-cycles by ith co-boundaries $H^i=Z^i/B^i$. Upper index of
an arbitrary form will denote its order, i.e., $\omega ^{k}\in
\Lambda ^{k}(\mathbb{R}^n)$. For standard de Rham complex in
$\mathbb{R}^n$ we take $\omega ^{k}\in B^k \subset \Lambda
^k(\mathbb{R}^n)$. Then $\omega ^{k}=d\nu ^{k-1}$ and
\[
H_k^k(\mathbb{R}^n)=Z^k/B^k.
\]
Let us explain the origin of the lower index in
$H_k^k(\mathbb{R}^n)$ \cite{5}. By definition,
\[Z^k=\ker \left( d:\Lambda
^{k}(\mathbb{R}^n)\rightarrow \Lambda ^{k+1}(\mathbb{R}^n)\right),
\qquad B^k=\text{im}\left( d:\Lambda ^{k-1}(\mathbb{R}^n)\rightarrow
\Lambda ^{k}(\mathbb{R}^n)\right).\] In $\Lambda
^{k}(\mathbb{R}^n)$, there exists forms
\[
\omega ^{k}=\bigwedge_{i=1}^{k}\lambda _{i}^{1},\text{\;  where \;
}\lambda ^{1}\in \Lambda ^{1}(\mathbb{R}^n)
\]
such, that $d\lambda ^{1}=0$. Suppose that $\lambda ^1 \in B^1$
(i.e., $\lambda ^1=d\mu ^0)$. Then
\[
d\omega ^k=0, \qquad \omega ^k=\bigwedge_{i=1}^{k}d\mu^0_i,
\]
i.e., exact form $\omega ^{k}\in \Lambda ^{k}(\mathbb{R}^n)$ is
itself an exterior product of exact forms. We obtain the quotient
space which we denote by
\[
H_{1,1,...,1}^{k}(\mathbb{R}^n)=Z^k/B^{1,1,...,1},
\]
where $B^{1,1,...,1}=\bigoplus%
\limits_{i=1}^{k}\text{im}_{i}\left( d:\Lambda
^{0}(\mathbb{R}^n)\rightarrow \Lambda ^{1}(\mathbb{R}^n)\right)$.

In similar way we get:

\noindent for $\omega ^{k}=d\mu ^1\wedge \bigwedge_{i=1}^{k-2}d\mu
_{i}^{0}$%
\[
H_{2,1,1,...,1}^{k}(\mathbb{R}^n)=Z^k/B^{2,1,1,...,1},
\]
where $B^{2,1,1,...,1}=\text{im}\left(
d:\Lambda ^{1}(\mathbb{R}^n)\rightarrow \Lambda ^{2}(\mathbb{R}^n)\right) \bigwedge\limits_{i=1}^{k-2}\text{im}_{i}\left( d:\Lambda ^{0}(%
\mathbb{R}^n)\rightarrow \Lambda ^{1}(\mathbb{R}^n)\right)$;

\noindent for $\omega ^{k}=d\mu ^{2}\wedge \bigwedge_{i=1}^{k-3}d\mu
_{i}^{0}$%
\[
H_{3,1,1,...,1}^{k}(\mathbb{R}^n)=Z^k/B^{3,1,1,...,1},
\]
where $B^{3,1,1,...,1}=\text{im}\left(
d:\Lambda ^{2}(\mathbb{R}^n)\rightarrow \Lambda ^{3}(\mathbb{R}^n)\right) \bigwedge\limits_{i=1}^{k-3}\text{im}_{i}\left( d:\Lambda ^{0}(%
\mathbb{R}^n)\rightarrow \Lambda ^{1}(\mathbb{R}^n)\right)$;

\noindent for $\omega ^{k}=d\mu ^{3}\wedge \bigwedge_{i=1}^{k-4}d\mu
_{i}^{0}$
\[
H_{4,1,1,...,1}^{k}(\mathbb{R}^n)=Z^k/B^{4,1,1,...,1},
\]
where $B^{4,1,1,...,1}=\text{im}\left(
d:\Lambda ^{3}(\mathbb{R}^n)\rightarrow \Lambda ^{4}(\mathbb{R}^n)\right) \bigwedge\limits_{i=1}^{k-4}\text{Im}_{i}\left( d:\Lambda ^{0}(%
\mathbb{R}^n)\rightarrow \Lambda ^{1}(\mathbb{R}^n)\right)$;

\begin{center}
...   ...   ...
\end{center}

\noindent for $\omega ^{k}=d\mu^1_1\wedge d\mu^1_2\wedge
\bigwedge_{i=1}^{k-4}d\mu _{i}^{0}$%
\[
H_{2,2,1,...,1}^{k}(\mathbb{R}^n)=Z^k/B^{2,2,1,...,1},
\]
where $B^{2,2,1,...,1}=%
\bigwedge\limits_{i=1}^{2}\text{Im}_{i}\left( d:\Lambda
^{1}(\mathbb{R}^n)\rightarrow \Lambda ^{2}(\mathbb{R}^n)\right)
\bigwedge\limits_{i=1}^{k-4}\text{Im}_{i}\left( d:\Lambda
^{0}(\mathbb{R}^n)\rightarrow \Lambda ^{1}(\mathbb{R}^n)\right)$;

\begin{center}
...   ...   ...
\end{center}
etc.

Cleary, for any Namby phase flow $\omega \in \Lambda
^{n-1}(\mathbb{R}^n)$, we have $[\omega] \in
H_{1,1,...,1}^{n-1}(\mathbb{R}^n)$. Moreover, by definition $\dim
H_{1,1,...,1}^{n-1}(\mathbb{R}^n)=0$.

\bigskip
\textbf{Theorem 2.} Let $\omega \in \Lambda ^{n-1}(\mathbb{R}^n)$ be
a divergence-free system. Then $[\omega] \in
H_{1,1,...,1}^{n-1}(\mathbb{R}^n)$.

\bigskip
\textbf{Proof.} Clearly, we need divergence-free condition only to
have $\omega$ being closed. Further proof is by induction.

For $\omega^{n-1} \in \Lambda^{n-1}(\mathbb{R}^n)$, we have
decomposition
\[
\omega^{n-1}=\sum\limits_{i=1}^n  A_idS_i^n
=\sum\limits_{i=1}^{n-1}A_idS_i^{n-1} \wedge
\left(dx_n-\frac{A_n}{A_1}\;dx_1\right) \qquad \text{ò.å.} \qquad
\omega^{n-1}=\omega^{n-2} \wedge \lambda^1,
\]
where $\omega^{n-2} \in \Lambda^{n-2}(\mathbb{R}^{n-1})$. Upper
index of hypersurface $S^{n}$ denotes the dimension of surround
space. Descedin in dimension till $\omega^2 \in \Lambda
^2(\mathbb{R}^3)$, we always have
\[
\sum\limits_{i=1}^3
A_idS_i^3=(A_1dx_2-A_2dx_1)\wedge\left(dx_3-\frac{A_3}{A_1}\;dx_1\right)
\qquad \text{ò.å.} \qquad \omega^2=\lambda^1_1 \wedge
\lambda^1_2.\quad \Box
\]

Theorem 2 provides the direct formula
\[
\sum\limits_{i=1}^n A_idS_i =A_1 \bigwedge \limits_{i=2}^{n}
\left(dx_i-\frac{A_i}{A_1}\;dx_1\right),\qquad \text{i.e.,} \qquad
\omega^{n-1}=\bigwedge \limits^{n-1}_{i=1} \lambda^1_i. \eqno{(4)}
\]
From the point of view of algebraic geometry, formula (4) can be
considered as one of the invariants of real coordinate realization
of of Birkhoff-Grothendieck theorem which states that on a
projective curve, all holomorphic vector bundles split into direct
sum of line bundles.

Further we will extesively use the analogy between (3)and (4).

\textbf{Definition.} Divergence-free phase flow
\[
\overset{\cdot}{\textbf{x}}=\{\textbf{h},\textbf{x}\}. \eqno{(5)}
\]
is called \textbf{Nambu mechanics}, if $[d\textbf{h}] \in
H_{1,1,...,1}^{n-1}(\mathbb{R}^n)$ wherein the class
$[d\textbf{h}]=0$; \textbf{Poincar\'e mechanics}, if $[d\textbf{h}]
\in H_{1,1,...,1}^{n-1}(\mathbb{R}^n)$ (wherein the class
$[d\textbf{h}] \neq 0$), or \textbf{Cartan mechanics} in remaining
of the cases.

\section{Non-integrable system}

Non-integrability of the system of equations
\[
\left\{\begin{array}{l}
\overset{\cdot}{x}= z^2  \\
\overset{\cdot}{y}= x^2 \\
\overset{\cdot}{z}= y^2
                      \end{array}\right.
\]
was discussed, e.g., in \cite{6}. This system can be rewritten in
the form (5) with vector Hamiltonian
\[
\textbf{h}=\frac{1}{4}\left(\begin{array}{c}
                        x^2z-y^3 \\
                        xy^2-z^3 \\
                        yz^2-x^3
                      \end{array}\right).
\]
Due to (3), (4)
\[
d\textbf{h}=\textbf{J}_1 \wedge \textbf{J}_2 \qquad \text{ãäå}
\qquad \textbf{J}_1=z^2dy-x^2dx, \qquad
\textbf{J}_2=dz-\frac{y^2}{z^2}\,dx,
\]
i.e., $[d\textbf{h}] \in H_{1,1}^{2}(\mathbb{R}^3)$ is Poincar\'e
mechanics.

\section{Symplectic mechanics}

The Hamiltonian
\[
H=\frac{1}{2}\sum\limits_{i=1}^2(x_i^2+p_i^2)
\]
generates integrable equations in $\mathbb{R}^4$
\[
\overset{\cdot}{x}_i=p_i, \qquad \overset{\cdot}{p}_i=-x_i.
\]
From Nambu point of view, to be completely integrable, this phase
flow should have three invariants
\[
d\textbf{h}=\textbf{J}_1 \wedge \textbf{J}_2\wedge \textbf{J}_3.
\]
The form of the first two integrals is due to exactness
$\textbf{J}_i=d\textbf{I}_i$, i.e.,
\[
\textbf{I}_i=\frac{1}{2}(x_i^2+p_i^2), \quad \text{äëÿ} \quad i=1,2.
\]
The third invariant has the form
\[
\textbf{J}_3=\frac{1}{2}\left(
\frac{dx_1}{p_1}-\frac{dp_1}{x_1}-\frac{dx_2}{p_2}+\frac{dp_2}{x_2}
\right).
\]
Then $d\textbf{h}=d\textbf{I}_1 \wedge d\textbf{I}_2\wedge
\textbf{J}_3$, where
\begin{eqnarray*}
\textbf{h}&=&\frac{1}{4}((p_2^2+x_2^2)\,dx_1\wedge
dp_1+(x_1^2+p_1^2)\,dx_2\wedge dp_2- (p_2p_1+x_1x_2)\,dx_2\wedge
dp_1  \\
&+&(x_1p_2-x_2p_1)\,dp_1\wedge dp_2-(p_2p_1+x_1x_2)\,dx_1\wedge
dp_2+ (x_1p_2-x_2p_1)\,dx_1\wedge dx_2).
\end{eqnarray*}
Since $d\textbf{J}\neq 0$, then $[d\textbf{h}] \in
H_{1,1,1}^3(\mathbb{R}^4)$ and $\dim H_{1,1,1}^3(\mathbb{R}^4) \neq
0$.

For symplectic mechanics in $\mathbb{R}^6$, we get $d\textbf{h}=
\bigwedge \limits_{i=1}^3 d\textbf{I}_i \wedge \textbf{J}$, where
\[ \textbf{J}=\sum\limits_{i>j}^3 \textbf{K}_i
\wedge \textbf{K}_j, \qquad \textbf{K}_i=\frac{1}{2}\left(
\frac{dx_i}{p_i}-\frac{dp_i}{x_i} \right), \quad i=1,2,3.
\]

In general case for symplectic mechanics in $\mathbb{R}^{2n}$,
 we get $d\textbf{h}= \bigwedge \limits_{i=1}^n d\textbf{I}_i \wedge
\textbf{J}$, where
\[
\textbf{J}=\sum\limits_i^n (-)^i\textbf{K}_1 \wedge ... \wedge
[\textbf{K}_i] \wedge ... \wedge \textbf{K}_n.
\]
Since $d\textbf{J}\neq 0$, then $[d\textbf{h}] \in
H_{n-1,1,...,1}^{2n-1}(\mathbb{R}^{2n})$, and $\dim
H_{n-1,1,...,1}^{2n-1}(\mathbb{R}^{2n}) \neq 0$, i.e., symplectic
mechanics in $\mathbb{R}^{2n}$ is the Cartan mechanics.

\section{Euler top}

The motion of a solid body around a fixed point is described by the
following equations in $\mathbb{R}^6$ \cite{7}:
\begin{eqnarray*}
\overset{\cdot}{x}_i&=&  \varepsilon_{ijk}\,\frac{x_j}{j_j}\, x_k+\varepsilon_{ijk}X_j\, y_k \\
\overset{\cdot}{y}_i&=& \varepsilon_{ijk} \,y_j\, x_k, \qquad
i,j,k=1,2,3.
\end{eqnarray*}
For all \textbf{X} and \textbf{j}, this system has three invariants:
\[
\textbf{I}_1=y_1^2+y_2^2+y_3^2, \qquad
\textbf{I}_2=x_1y_1+x_2y_2+x_3y_3, \qquad
\textbf{I}_3=\frac{x_1^2}{j_1}+\frac{x_2^2}{j_2}+\frac{x_3^2}{j_3}+X_1y_1+X_2y_2+X_3y_3.
\]
In the case of $X_1=X_2=X_3=0$, we get another two invariants:
\[
\textbf{I}_4=x_1^2+x_2^2+x_3^2, \qquad
\textbf{J}_5=\frac{1}{\Delta}\left(
\frac{x_1dy_1}{j_1}+\frac{x_2dy_2}{j_2}+\frac{x_3dy_3}{j_3} \right),
\]
where $\Delta= \overset{\cdot}{\textbf{x}} \cdot
\textbf{y}=\varepsilon_{ijk}\,y_i\frac{x_j}{j_j}\, x_k$. Thus
\[
d\textbf{h}= \bigwedge \limits_{i=1}^4 d\textbf{I}_i \wedge
\textbf{J}_5,
\]
i.e., $[dh] \in H_{1,1,1,1,1}^{5}(\mathbb{R}^6)$, and Euler top is
Poincar\'e mechanics.

\bigskip
\bigskip
\begin{flushright}
{\footnotesize Translated by P.I.Troshin}
\end{flushright}

\end{document}